\documentclass[12pt]{article}
\usepackage{amsfonts}
\usepackage{amsmath}

\newtheorem{theorem}{Theorem}

\newtheorem{propo}{Proposition}

\newenvironment{proof}[1][Proof]{\noindent\textbf{#1.} }{\ \rule{0.5em}{0.5em}}

\def\R{\mathbb R_q^+}

\begin{document}

\title{\bf Applications of the $q$-Fourier Analysis to the Symmetric Moment Problem}
\date{ }
\author{Lazhar Dhaouadi \thanks{%
Institut Pr\'eparatoire aux Etudes d'Ing\'enieur de Bizerte
(Universit\'e du 7 novembre à Carthage). Route Menzel Abderrahmene
Bizerte, 7021 Zarzouna, Tunisia. \quad\quad\quad\quad\quad\quad
E-mail lazhardhaouadi@yahoo.fr}} \maketitle

\begin{abstract}
Sufficient condition for the symmetric moment problem to be
determinate is given using standards methods of $q$-Fourier
analysis. This condition it cannot be a particular case of
Carleman's criterion.
\end{abstract}

\section{Introduction}

For a positive measure $\mu$ on $\mathbb R$, the nth moment is
defined as
$$
s_n=\int x^nd\mu(x),
$$
provided the integral exists. Let $\mathcal M^*(\mathbb R)$ be the
set of positive measures on $\mathbb R$ with moment of any order and
with infinite support. Given $\mu\in\mathcal M^*(\mathbb R)$, then
we consider $V_\mu$ the set of all $\nu\in\mathcal M^*(\mathbb R)$
such that
$$
s_n=\int x^nd\mu(x)=\int x^nd\nu(x)\quad\textrm{for}\quad n\geq 0.
$$
We say that $\mu$ is determinate if $V_\mu=\{\mu\}$, otherwise $\mu$
is indeterminate. For more detail to this subject, the reader can
consult the references [1] and [2], and the references theine. We
recall that a moment problem is said to be symmetric if all moments
of odd order are $0$. It is desirable to be able to tell whether the
moment problem is determinate or indeterminate just by looking at
the moment sequence $(s_n)_{n\geq0}$. We have the following
classical sufficient condition for determinacy:

\bigskip
The Perron's criterion
\begin{equation}
\limsup_{n\rightarrow\infty}\left\{\frac{s_{2n}}{(2n)!}\right\}^{1/2n}<\infty
\end{equation}

The Riesz's criterion
\begin{equation}
\liminf_{n\rightarrow\infty}\left\{\frac{s_{n}}{n!}\right\}^{1/n}<\infty
\end{equation}

The Carleman's criterion
\begin{equation}
\sum_{n=1}^\infty\frac{1}{\sqrt[2n]{s_{2n}}}=\infty
\end{equation}
It follows easily that we have
$$
(1)\Rightarrow (2)\Rightarrow (3).
$$
So by Carleman's theorem (see [3]), if $\mu$ satisfy (3) then $\mu$
is determinate. In the end we have the following equivalent form of
the Perron condition
$$
\int_0^\infty e^{\alpha|x|}d\mu(x)<\infty,\quad\text{for some}\quad
\alpha>0.
$$
For $\mu\in\mathcal M^*(\mathbb R)$ let $(P_n)$ be the corresponding
orthonormal polynomials determined by
$$
P_n(x)=k_nx^n+....,\quad k_n>0;\quad \int
P_n(x)P_m(x)d\mu(x)=\delta_{nm}.
$$
For an indeterminate $\mu\in\mathcal M^*(\mathbb R)$, it is known
that the series $\sum|P_n(z)|^2$ and $\sum|Q_n(z)|^2$ converge
uniformly on compact subsets of $\mathbb C$, where $(Q_n)$ the
polynomials of the second kind
$$
Q_n(x)=\int\frac{P_n(x)-P_n(y)}{x-y}d\mu(y),\quad x\in\mathbb C.
$$
Therefore, the series
$$\aligned
A(z)&=z\sum_{n=0}^\infty Q_n(0)Q_n(z)\\
B(z)&=-1+z\sum_{n=0}^\infty Q_n(0)P_n(z)
\\
C(z)&=1+z\sum_{n=0}^\infty P_n(0)Q_n(z)
\\
D(z)&=z\sum_{n=0}^\infty P_n(0)P_n(z)
\endaligned$$
determine entire functions. For $t\in\mathbb R\cup\{\infty\}$, we
introduce the special discrete measures so called N-extremal of the
form
$$\aligned
\nu_t&=\sum_{\lambda\in\Lambda_t}m_\lambda\delta_\lambda,\quad\Lambda_t=\{z\in\mathbb
C\big|B(z)t-D(z)=0\},\\
m_\lambda&=\frac{A(\lambda)t-C(\lambda)}{B'(\lambda)t-D'(\lambda)}=\left(\sum_{n=0}^\infty
P_n^2(\lambda)\right)^{-1},\quad \text{for}\quad\lambda\in\Lambda_t,
\endaligned$$
with $\Lambda_\infty=\{z\in\mathbb C\big|B(z)=0\}$.

\bigskip
We recall the following result of M. Riesz [6]:

\bigskip
Let $\mu\in\mathcal M^*(\mathbb R)$
\begin{enumerate}
\item If $\mu$ is indeterminate and $\nu\in V_\mu$ then $(P_n)$
form an orthonormal basis of $L^2(\nu)$ if and only if $\nu$ is
N-extremal.
\item If $\mu$ is determinate then $(P_n)$
form an orthonormal basis of $L^2(\mu)$.
\end{enumerate}

\section{Main Results}

\bigskip
We consider $0<q<1$ and we design by $\mu$  the positive measure
defined by
$$
d\mu(x)=\omega(x)^2d\sigma(x)=\omega(x)^2|x|^{2v+1}d_qx,\quad v>-1,
$$
where $d_qx$ a discrete measure defined by
$$
d_qx=(1-q)\sum_{n\in\mathbb Z}q^n\delta_{q^n}.
$$
and $\omega$ an even real function satisfies $\omega(q^n)\neq 0$ for
all $n\in\mathbb Z$.

\bigskip
Let $(P_n)$ be the corresponding orthonormal polynomials for $\mu$,
and we consider $\R$ the support of $\mu$
$$
\R=\{q^n,\quad n\in\mathbb Z\}.
$$

\begin{theorem}
The sequence $(P_n)$ form an orthonormal basis of $L^2(\mu)$ if and
only if $\mu$ is determinate.
\end{theorem}

\begin{proof} If the sequence $(P_n)$ form an orthonormal basis of
$L^2(\mu)$ then $\mu$ is determinate or N-extremal. Suppose that
$\mu$ is N-extremal, then there exist $t\in\mathbb R\cup\{\infty\}$
such that
$$
\mu=\nu_t=\sum_{\lambda\in\Lambda_t}m_\lambda\delta_\lambda.
$$
This implies $\Lambda_t=\R$. In the other hand $\Lambda_t$ is the
set of zero of entire function and $\{0\}$ is an accumulation point
of $\Lambda_t$, then $\Lambda_t=\mathbb C$, which is absurd.
\end{proof}

\begin{propo}
The sequence $(P_n)$ form an orthonormal basis of $L^2(\mu)$ if and
only if
$$
\psi_n(x)=\omega(x)P_n(x),
$$
form an orthonormal basis of $L^2(\sigma)$.
\end{propo}

\begin{proof} It is clear that
$$
\int P_n(x)P_m(x)d\mu(x)=\delta_{nm}\Leftrightarrow\int
\psi_n(x)\psi_m(x)d\sigma(x)=\delta_{nm}.
$$
Now, let $f$ be a function belongs to $L^2(\sigma)$ and  put
$g(x)=\frac{f(x)}{\omega(x)}$, then $g\in L^2(\mu)$ and we have
$$
\int g(x)P_n(x)d\mu(x)=0\Leftrightarrow \int
f(x)\psi_n(x)d\sigma(x)=0,
$$
which lead to the result.
\end{proof}

\begin{theorem}Let $s_n$ be the nth moment of $\sigma$. If
\begin{equation}
\lim_{n\rightarrow\infty}q^{n/4}\sqrt[2n]{s_{2n}}=0,
\end{equation}
then $\mu$ is determinate.
\end{theorem}

\begin{proof}Using theorem 2, it suffice to prove that $(P_n)$ form
an orthonormal basis of $L^2(\mu)$ and by the proposition 1, this is
equivalent to prove that $(\psi_n)$ form an orthonormal basis of
$L^2(\sigma)$. Let $f$ be a function belongs to $L^2(\sigma)$ such
that
$$
\int f(x)\psi_n(x)d\sigma(x)=0,\quad\forall n\in\mathbb N.
$$
We can write $f=f^++f^-$, where $f^+$ is even and $f^-$ is odd. This
implies
$$
\int f^+(x)x^{2n}\omega(x)d\sigma(x)=\int
f^-(x)x^{2n+1}\omega(x)d\sigma(x)=0,\quad\forall n\in\mathbb N.
$$
Then
$$
\int f^+(x)\omega(x)j_v(\lambda x,q^2)d\sigma(x)=\int
xf^-(x)\omega(x)j_v(\lambda x,q^2)d\sigma(x)=0,\quad\forall
\lambda\in\mathbb R_q^+,
$$
where $j_v(x,q^2)$ is the the normalized $q$-Bessel function of
third kind defined by
$$
j_v(x,q^2)=\sum_{n=0}^\infty (-1)^n
\frac{q^{n(n+1)}}{(q^2,q^2)_n(q^{2v+2},q^2)_n}x^{2n}.
$$
We can exchange integral and sum because
$$\aligned
&\sum_{n=0}^\infty
\frac{q^{n(n+1)}}{(q^2,q^2)_n(q^{2v+2},q^2)_n}\lambda^{2n}\int
|f^+(x)|x^{2n}\omega(x)d\sigma(x)\\
&\leq\|f^+\|_2\sum_{n=0}^\infty
\frac{q^{n(n+1)}}{(q^2,q^2)_n(q^{2v+2},q^2)_n}\lambda^{2n}\left[\int
x^{4n}\omega^2(x)d\sigma(x)\right]^{1/2}\\
&\leq\|f^+\|_2\sum_{n=0}^\infty
\frac{q^{n(n+1)}}{(q^2,q^2)_n(q^{2v+2},q^2)_n}\sqrt{s_{4n}}\lambda^{2n}<\infty.
\endaligned$$
and
$$\aligned
&\sum_{n=0}^\infty
\frac{q^{n(n+1)}}{(q^2,q^2)_n(q^{2v+2},q^2)_n}\lambda^{2n}\int
|xf^-(x)|x^{2n}\omega(x)d\sigma(x)\\
&\leq\|f^-\|_2\sum_{n=0}^\infty
\frac{q^{n(n+1)}}{(q^2,q^2)_n(q^{2v+2},q^2)_n}\lambda^{2n}\left[\int
x^{2(2n+1)}\omega^2(x)d\sigma(x)\right]^{1/2}\\
&\leq\|f^-\|_2\sum_{n=0}^\infty
\frac{q^{n(n+1)}}{(q^2,q^2)_n(q^{2v+2},q^2)_n}\sqrt{s_{2(2n+1)}}\lambda^{2n}<\infty.
\endaligned$$
The convergence of the above series is a simple consequence of the
Cauchy root test (4).

\bigskip
In [4] the authors introduce the $q$-Bessel Fourier transform as
follows
$$
\mathcal F_{q,v}f(\lambda)=c_{q,v}\int_0^\infty f(x)j_v(\lambda
x,q^2)x^{2v+1}d_qx,\quad\forall \lambda\in\mathbb R_q,
$$
where $f$ is an even function belongs to $L^1(\sigma)$ and he have
proved that if $\mathcal F_{q,v}f=0$ then $f=0$. Using the
Cauchy-Schwartz inequality we prove that $f^+\omega$ and
$xf^-\omega$ belongs to $L^1(\sigma)$, which implies that
$f^+=f^-=0$.
\end{proof}

\section{Application}
In [5], the Ramanujan identity was proved
$$
\sum_{k\in\mathbb
Z}\frac{z^{k}}{(bq^{k},q)_{\infty}}=\frac{(bz,q/bz,q,q)_{\infty}
}{(b,z,q/b,q)_{\infty}},
$$
where
$$
(x,q)_{\infty}=\prod_{i=0}^{\infty}(1-q^ix).
$$
Let $\alpha>0$ and replace in the above sum
$$
b\rightarrow-1,z\rightarrow q^{\alpha}z,
$$
we obtain
$$
\sum_{k\in\mathbb Z} q^{\alpha
k}\frac{z^{k}}{(-q^{k},q)_{\infty}}=\frac{(-q^{\alpha
}z,-q^{-\alpha}q/z,q,q)_{\infty}}{(-1,q^{\alpha}z,-q,q)_{\infty}}%
=\frac{(-q^{\alpha}z,q,q)_{\infty}}{(-1,q^{\alpha}z,-q,q)_{\infty}%
}(-q^{-\alpha}q/z,q)_{\infty}.
$$
On the other hand
\begin{align*}
(-q^{-\alpha}q/z,q)_{\infty} &  =\prod_{i=0}^{\infty}(1+q^{i-\alpha}%
q/z)=\prod_{i=0}^{[\alpha]}(1+q^{i-\alpha}q/z)\prod_{i=[\alpha]+1}^{\infty
}(1+q^{i-\alpha}q/z)\\
&  =\left(  \frac{q}{z}\right)
^{\alpha}\prod_{i=0}^{[\alpha]}q^{i-\alpha
}\prod_{i=0}^{[\alpha]}(1+q^{\alpha-i}z/q)\prod_{i=[\alpha]+1}^{\infty
}(1+q^{i-\alpha}z/q)\\
&  =\left(  \frac{q}{z}\right)  ^{\alpha}\prod_{i=0}^{[\alpha]}(1+q^{i}%
q^{\alpha-[\alpha]}z/q)\prod_{i=[\alpha]+1}^{\infty}(1+q^{i-\alpha
}z/q)q^{\left(  \frac{[\alpha]}{2}-\alpha\right)  \left(
[\alpha]+1\right)  }.
\end{align*}
This gives
\begin{align*}
& \sum_{k\in\mathbb Z} q^{\alpha k}\frac{z^{k}}{(-q^{k},q)_{\infty}}\\
& =\left[  \frac{(-q^{\alpha}z,q,q)_{\infty}}{(-1,q^{\alpha}z,-q,q)_{\infty}%
}\prod_{i=0}^{[\alpha]}(1+q^{i}q^{\alpha-[\alpha]}z/q)\prod_{i=[\alpha
]+1}^{\infty}(1+q^{i-\alpha}z/q)\right]  \left(  \frac{q}{z}\right)
^{\alpha }q^{\left(  \frac{[\alpha]}{2}-\alpha\right)  \left(
[\alpha]+1\right)  }.
\end{align*}
Put
$$
C(\alpha,q,z)=\left[
\frac{(-q^{\alpha}z,q,q)_{\infty}}{(-1,q^{\alpha
}z,-q,q)_{\infty}}\prod_{i=0}^{[\alpha]}(1+q^{i}q^{\alpha-[\alpha]}%
z/q)\prod_{i=0}^{\infty}(1+q^{i+[\alpha]-\alpha}z)\right]  \left(  \frac{q}%
{z}\right)  ^{\alpha},
$$
and
$$
\sigma(\alpha)=\left(  \frac{[\alpha]}{2}-\alpha\right)  \left(
[\alpha]+1\right)
$$
Now, let $p>0$ and replace
$$
\alpha\rightarrow\frac{n}{p},q\rightarrow q^{2p},z\rightarrow
q^{2v+2},
$$
we get
$$
(1-q)\sum_{k\in\mathbb Z}
q^{2nk}\frac{q^{2k(v+1)}}{(-q^{2pk},q^{2p})_{\infty}}=(1-q)C\left(
\frac{n}{p},q^{2p},q^{2v+2}\right)  q^{2p\sigma(\frac{n}{p})},
$$
We recall the definition of the $q$-exponential function
$$
e(x,q)=\frac{1}{(x,q)_{\infty}}=\prod_{i=0}^{\infty}\frac{1}{1-q^ix},
$$
and we can write
$$
s_{2n}=\int_{0}^{\infty}x^{2n}e(-x^{2p},q^{2p})x^{2v+1}d_{q}x=(1-q)C\left(
\frac{n}{p},q^{2p},q^{2v+2}\right)  q^{2p\sigma(\frac{n}{p})},
$$
and then
$$
\sqrt[2n]{s_{2n}}=\left[  (1-q)C\left(
\frac{n}{p},q^{2p},q^{2v+1}\right) \right]
^{\frac{1}{2n}}q^{\frac{p}{n}\sigma(\frac{n}{p})}.
$$
It is easy to prove that
$$
\lim_{n\rightarrow\infty}\left[  (1-q)C\left(  \frac{n}{p},q^{2p}%
,q^{2v+2}\right)  \right]  ^{\frac{1}{2n}}=1.
$$
For the particular case
$$
p=3,
$$
we have
$$
\lim_{n\rightarrow\infty}\frac{3}{n}\sigma\left(  \frac{n}{3}\right)
+\frac{n}{4}=+\infty\Rightarrow\lim_{n\rightarrow\infty}q^{n/4}\sqrt[2n]%
{s_{2n}}=0
$$
and there exist $n_0>0$ such that for all $n>n_0$
$$
\frac{3}{n}\sigma\left(  \frac{n}{3}\right)
+\frac{n}{4}>0\Rightarrow \sum_{n=0}^\infty
\frac{1}{\sqrt[2n]{s_{2n}}}<\infty.
$$
In the end, if we denote by
$$
d\mu(x)=e(-x^6,q^6)x^{2v+1}d_qx,
$$
then by the Theorem 2 we see that $\mu$ is determinate, but not
satisfying the Carleman's criterion.

\end{document}